# Introduction to Projective Arithmetics

## M. Burgin

Department of Mathematics
*University of California, Los Angeles*
Los Angeles, CA 90095

**Abstract:**

Science and mathematics help people to better understand world, eliminating many inconsistencies, fallacies and misconceptions. One of such misconceptions is related to arithmetic of natural numbers, which is extremely important both for science and everyday life. People think their counting is governed by the rules of the conventional arithmetic and thus other kinds of arithmetics of natural numbers do not exist and cannot exist. However, this popular image of the situation with the natural numbers is wrong. In many situations, people have to utilize and do implicitly utilize rules of counting and operating different from rules and operations in the conventional arithmetic. This is a consequence of the existing diversity in nature and society. To correctly represent this diversity, people have to explicitly employ different arithmetics. To make a distinction, we call the conventional arithmetic by the name *Diophantine arithmetic*, while other arithmetics are called non-Diophantine. There are two big families of *non-Diophantine arithmetics*: *projective arithmetics* and *dual arithmetics* (Burgin, 1997). In this work, we give an exposition of projective arithmetics, presenting their properties and considering also a more general mathematical structure called a projective prearithmetic. The Diophantine arithmetic is a member of this parametric family: its parameter is equal to the identity function $f(x) = x$. In conclusion, it is demonstrated how non-Diophantine arithmetics may be utilized beyond mathematics and how they allow one to eliminate inconsistencies and contradictions encountered by other researchers.

## 1. Introduction

From all things that people are doing, counting is one of the most important. Without counting, people cannot do a lot of things: people cannot buy and sell, people cannot develop science and technology, people cannot organize mass production and so on and so forth. Every woman and every man, every boy and every girl perform counting many times a day. Calculators and computers were invented to help people to count. Later computer began to fulfill much more sophisticated tasks. Some of its abilities look miraculous. However, counting lies at the bottom of all computer operations.

People started counting using numbers in prehistoric times. Numbers used for counting are called natural numbers or counting numbers. These numbers are the most basic for human civilization. All other numbers, such as whole numbers, rational numbers, real numbers, complex numbers, etc., were invented much later than natural numbers. The system of natural numbers with operations is called arithmetic of natural numbers. However, because this is the most profound arithmetic, we call it simply arithmetic or the Diophantine arithmetic. Although there are arithmetics of whole numbers, rational numbers, real numbers, complex numbers, and so on, operations with natural numbers induce operations in all other arithmetics.

For thousands of years mathematicians studied natural numbers and counting learning a lot in this area. People's experience with numbers has been profound but incomplete.

The situation with arithmetic now is similar to the situation with geometry in the middle of the $19^{th}$ century. Namely, the Euclidean geometry was believed for 2200 years to be unique (both as an absolute truth and a necessary mode of human perception). People were not even able to imagine something different. The famous German philosopher Emmanuil Kant claimed that (Euclidean) geometry is given to people *a priory*, i.e., without special learning. In spite of this, almost unexpectedly some people began to understand that geometry is not unique. Trying to improve the axiomatic system suggested for geometry by Euclid, three great mathematicians of the 19th century (C.F.Gauss, N.I. Lobachewsky, and Ja. Bolyai) discovered a lot of other, non-Euclidean geometries. At first, even the best mathematicians opposed this discovery and severely attacked Lobachewsky and Bolyai who published their results. Forecasting such antagonistic attitude, the first mathematician of his times Gauss was afraid to publish his results on non-Euclidean geometry. Nevertheless, progress of mathematics brought understanding and then recognition. This discovery is now considered as one of the highest achievements of the human genius. It changed to a great extent understanding of mathematics and improved comprehension of the whole world.

The situation with arithmetic is even more striking. For thousands of years, much longer than that for the Euclidean geometry, people explicitly used only one arithmetic and thought that nothing else could exist. Mathematical establishment treated arithmetic as primordial entity. For example, such prominent German mathematician as Leopold Kronecker (1825-1891) wrote: "*God made the integers, all the rest is the work of man*". Carl Friedrich Gauss (1777-1855) claimed that " God does arithmetic."

Almost all people, mathematicians, as well as non-mathematicians, have had and have no doubts that 2 + 2 = 4 is the most evident truth in the world, which is valid always and everywhere. In the authoritative mathematical journal "The American Mathematical Monthly" (April, 1999, p.375), "*Although other sciences and philosophical theories change their 'facts' frequently*, 2 + 2 *remains* 4."

However, sages of the ancient Greece started to doubt convenience of the conventional arithmetic. There was a group of philosophers, who were called Sophists and lived from the second half of the fifth century B.C.E. to the first half of the fourth century B.C.E.. Sophists asserted relativity of human knowledge and suggested many paradoxes, explicating complexity and diversity of real world.

The famous philosopher Zeno of Elea (490-430 B.C.E.), who was said to be a self-taught country boy, invented very impressive paradoxes, in which he questioned the popular knowledge and intuition related to such fundamental essences as time, space, and number.

An example of such reasoning is the *paradox of the heap* (or the *Sorites paradox* where *σωρος* is the Greek word for "heap"). It is possible to formulate this paradox in the following way.

1. One million grains of sand make a heap.

2. If one grain of sand is added to this heap, the heap stays the same.

3. However, when we add 1 to any natural number, we always get a new number.

As we know, Greek sages posed questions, but in many cases, including arithmetic, suggested no answers. As a result, for more than two thousand years these problems were forgotten and everybody was satisfied with the conventional arithmetic. In spite of all problems and paradoxes, this arithmetic has remained very and very useful.

In modern times, scientists and mathematicians returned to problems with arithmetic. The famous German researcher Herman Ludwig Ferdinand von Helmholtz (1821-1894) was may be the first scientist who questioned adequacy of the conventional arithmetic. In his "Counting and Measuring" (1887), Helmholtz considered an important problem of applicability of arithmetic, at that time it was only one arithmetic, to physical phenomena. This was a natural approach of a scientist, who even mathematics judged by the main criterion of science – observation and experiment.

His first observation was that as the concept of number is derived from some practice, usual arithmetic is applicable to these experiences. However, it is easy to find many situations when this is not true. To mention but a few described by Helmholtz, one raindrop added to another raindrop does not make two raindrops.

In a similar way, when one mixes two equal volumes of water, one at 40° Fahrenheit and the other at 50°, one does not get two volumes at 90°. Alike the conventional arithmetic fails to describe correctly the result of combining gases or liquids by volumes. For example (Kline, 1980), one quart of alcohol and one quart of water yield about 1.8 quarts of vodka.

Later the famous French mathematician Henri Lebesgue facetiously pointed out (cf. (Kline, 1980)) that if one puts a lion and a rabbit in a cage, one will not find two animals in the cage later on.

However, very few (if any) paid attention to the work of Helmholtz on arithmetic, and as still no alternative to the conventional arithmetic has been suggested, these problems were once more forgotten. Only much later, in the second part of the 20$^{th}$ century, mathematicians began to doubt once more the absolute character of the ordinary arithmetic, where 2 + 2 = 4 and two times two is equal to four. Scientists and mathematicians draw attention of the scientific community to the foundational problems of natural numbers and the conventional arithmetic. The most extreme assertion that there is only a finite quantity of natural numbers was suggested by Yesenin-Volpin (1960), who developed a mathematical direction called ultraintuitionism and took this assertion as one of the central postulates of ultraintuitionism. Other authors also considered arithmetics with a finite number of numbers, claiming that these arithmetics are inconsistent (cf., for example, (Meyer and Mortensen, 1984; Van Bendegem, 1994; Priest, 1997; 2000; Rosinger, 2008)).

Van Danzig had similar ideas but expressed in different way. In his article (1956), he argued that only some of natural numbers may be considered finite. Consequently, all other mathematical entities that are called traditionally natural numbers are only some expressions but not numbers. These arguments are supported and extended by Blehman, *et al* (1983).

Other authors are more moderate in their criticism of the conventional arithmetic. They write that not all natural numbers are similar in contrast to the presupposition of the conventional arithmetic that the set of natural numbers is uniform (Kolmogorov, 1961; Littlewood, 1953; Birkhoff and Barti, 1970; Rashevsky, 1973; Dummett 1975; Knuth, 1976). Different types of natural numbers have been introduced, but without changing the conventional arithmetic. For example, Kolmogorov (1961) suggested that in solving practical problems it is worth to separate *small, medium, large,* and *super-large* numbers. Rashevsky (1973) explicitly formulated the problem of constructing a new arithmetic of natural numbers.

Such new arithmetics were discovered in 1975 and their theory was published in (Burgin, 1977). The conventional arithmetic may be called the Diophantine arithmetic because the ancient Greek mathematician Diophantus was the first who made an essential contribution to arithmetic. He is best known for his *Arithmetica*, a work with an essential contribution to the solution of algebraic equations and to arithmetic. However, essentially nothing is known of his life and there has been much debate regarding the date at which he lived (sometimes between 150 BCE and 350 CE) although it is generally assumed that he did his work in the great city of Alexandria..

Thus, it is natural to call new arithmetics non-Diophantine. Different properties of *non-Diophantine arithmetics* and some applications are described in (Burgin, 1992; 1997; 1998; 2007). Various examples of utilization of non-Diophantine arithmetics in science and in everyday life are demonstrated in (Burgin, 2001).

Like non-Euclidean geometries of Lobachewsky, non-Diophantine arithmetics depend on a special parameter, although this parameter is not a number as in the case of Lobachewsky geometries. Arithmetics have a functional parameter, implying that properties of and laws of operations in non-Diophantine arithmetics depend on a definite function $f(x)$.

There are two big families of non-Diophantine arithmetics: projective and dual arithmetics (Burgin, 1997). Here we give an exposition of projective arithmetics, presenting their properties and considering also a more general mathematical structure called a *projective prearithmetic*. The Diophantine arithmetic is a member of the parametric family of projective arithmetics: its parameter is equal to the identity function $f(x) = x$. One of the interesting properties of projective arithmetics is that they allow us to formalize and make rigorous such concepts as *much smaller* (<<) and *much larger* (>>), which are frequently used in theoretical physics and other areas of applied mathematics on the intuitive level without a proper mathematical formalization. In conclusion, it is demonstrated how non-Diophantine arithmetics may be utilized beyond mathematics and how they allow one to eliminate inconsistencies and contradictions encountered by other researchers.

**Denotations**

Let **N** be the arithmetic of all natural numbers with 0 while +, · and ≤ be operations of addition and multiplication**,** and the natural order relation in **N,** correspondingly. The set of all natural numbers is denoted by $N$. The set of all real numbers is denoted by $R$. If $X$ is a set, then $1_X : X \to X$ is the identity mapping, i.e., $1_X(x) = x$. If $n$ is a natural (whole) number, then $Sn$ denotes the next number $n + 1$.

## 2. Projective Arithmetics

At first, we define prearithmetics. Informally, a *prearithmetic* is a subset of all natural numbers with, at least, two operations + (addition) and ∘ (multiplication), which are defined for all its elements. In what follows, an arbitrary prearithmetic is denoted by $\boldsymbol{A} = (A; +, \circ)$ where $A$ is a subset of $N$ and is called the set of the elements or numbers of $\boldsymbol{A}$. An *arithmetic* is a prearithmetic such that satisfies extra conditions. One of such conditions essential for applications is the linearity of the order of the numbers (i.e., elements) of the arithmetic. As an example of application where this condition is important, it is possible to take utility function used for decision making. However, non-Diophantine arithmetics are not defined axiomatically (by this or any other conditions). They are constructed by the procedure that is described below, using the conventional natural numbers.

In what follows, an arbitrary arithmetic is denoted by $\boldsymbol{A} = (A; +, \circ, \leq)$. The above condition means that all numbers in $\boldsymbol{A}$ have to be linearly ordered by the relation $\leq$.

Let us take two prearithmetics $\boldsymbol{A}_1 = (A_1; +_1, \circ_1)$ and $\boldsymbol{A}_2 = (A_2; +_2, \circ_2)$ as well as two functions $g: A_1 \rightarrow A_2$ and $h: A_1 \rightarrow A_2$.

**Definition 1.** A *prearithmetic* $\boldsymbol{A}_1 = (A_1; +_1, \circ_1)$ is called weakly projective (projective) with respect to a prearithmetic $\boldsymbol{A}_2 = (A_2; +_2, \circ_2)$ if ($hg = 1_{A1}$ and $hg = 1_{A2}$) the two operations in $\boldsymbol{A}_1$ are defined for any two numbers $a$ and $b$ as follows:

$$a +_1 b = h(g(a) +_2 g(b));$$
$$a \circ_1 b = h(g(a) \circ_2 g(b)).$$

**Definition 2.** The function $g$ is called the *projector* and the function $h$ is called the *coprojector* for the pair $(\boldsymbol{A}_1, \boldsymbol{A}_2)$.

If we have a set $B$, a prearithmetic $\boldsymbol{A} = (A; +, \circ)$ and two functions $g: B \rightarrow A$ and $h: A \rightarrow B$, then it is possible to define on $B$ the structure of the prearithmetic that is weakly projective with respect to $A$. It will be a unique prearithmetic with the projector $g$ and the coprojector $h$.

**Theorem 1.** If a prearithmetic $A_1 = (A_1 ; +_1 , \circ_1)$ is weakly projective with respect to a prearithmetic $A_2 = (A_2 ; +_2 , \circ_2)$ and $gh = 1_{A2}$, then the prearithmetic $A_2$ is projective with respect to a prearithmetic $A_1$.

Let us consider an arbitrary non-decreasing function $f: U \to R$ defined on a subset $U$ of $R$. It is possible to define the following two functions $f_T: M \to N$ and $f^T: P \to M$ where $M$, $P$ are some subsets of $N$.

$$f_T(x) = \inf \{n; n \in N \text{ and } n \geq f(x)\} = ]f(x)[ = \lceil f(x) \rceil ;$$

$$f^T(x) = \sup \{m; m \in N \text{ and } m \leq f^{-1}(x)\} = [f^{-1}(x)] = \lfloor f^{-1}(x) \rfloor .$$

**Definition 3.** A prearithmetic $A = (M ; + , \circ)$ is called a *projective prearithmetic* if it is weakly projective with respect to the conventional (Diophantine) arithmetic $\mathbf{N} = (N ; + , \cdot , \leq)$ with the *projector* $f_T(x)$ and the *coprojector* $f^T(x)$. The function $f(x)$ is called the *generator* of the *projector* $f_T(x)$, *coprojector* $f^T(x)$ and the prearithmetic $A$.

Operations in $A_1$ are defined for any two numbers $a$ and $b$ as follows:

$$a + b = f^T(f_T(a) + f_T(b)) ;$$

$$a \circ b = f^T(f_T(a) \cdot f_T(b)) .$$

Elements of $M$ are called numbers of $A$ and are denoted with subscript $\mu$, i.e., 2 in $A$ is denoted by $2_\mu$ and 5 in $A$ is denoted by $5_\mu$. Numbers of $A$ are ordered by the same order relation $\leq$ as they are ordered in $\mathbf{N}$.

**Example 1.** Let us take the set $M = \{1, 2, \ldots , m\}$ and define the projector $f_T: M \to N$ as the natural inclusion of $M$ into $N$. The coprojector $f^T: N \to N$ is defined by the formula $f^T(n) = k$ where $1 \leq k < m$ and $n - k$ is divisible by $m$. Then the projective prearithmetic $A = (M ; + , \circ)$ defined by these functions is the arithmetic of residues modulo $m$.

Let $A = (M; +, \circ)$ be a projective prearithmetic.

**Theorem 2.** If $f_T(0_\mu) = 0$, then for any $a_\mu \in A$ the equality $0_\mu + a_\mu = a_\mu$ is true if and only if $f(x)$ is a strictly increasing function.

This theorem gives necessary and sufficient conditions for $0_\mu$ to be a neutral element (zero) with respect to addition in the projective prearithmetic $A$.

**Proposition 1.** If $f_T(0_\mu) = 0$, then for any $a_\mu \in A$ the equality $0_\mu \circ a_\mu = a_\mu$ is true.

It means that the natural equality $f_T(0_\mu) = 0$ always imply that $0_\mu$ is a absorbing element (zero) with respect to multiplication in the projective prearithmetic $A$.

It implies the following definition.

**Definition 4.** A projective prearithmetic $A = (N; +, \circ, \leq)$ with the projector $f_T(x)$ and the coprojector $f^T(x)$ is called a *projective arithmetic* if the following conditions are satisfied:

1) $f_T(0_\mu) = 0$;

2) $f(x)$ is a strictly increasing function;

3) for any elements $a$ and $b$ from $U$ from $a \leq b$, we have $f_T(Sa) - f_T(a) \leq f_T(Sb) - f_T(b)$.

Here, $Sa$ is the element of $A$ that strictly follows $a$.

Thus, operations in $A$ are defined for any two its numbers $a_\mu$ and $b_\mu$ as follows:

$$a_\mu + b_\mu = f^T(f_T(a) + f_T(b));$$
$$a_\mu \circ b_\mu = f^T(f_T(a) \cdot f_T(b)).$$

Besides, a natural order relation is defined on $A$:

$$a_\mu \leq b_\mu \text{ if and only if } a \leq b.$$

**Remark 1.** It is possible to define the generator (function $f$) only for natural numbers, but in many cases, its definition for real numbers makes the analytic expression for $f$ simpler.

**Example 2.** Let us take as a generator $f(x)$ for a projective arithmetic $A = (N; +, \circ, \leq)$ a simple function, such as $x^2$, and look how operations are performed.

$2_\mu + 2_\mu = f^T(f_T(2_\mu) + f_T(2_\mu)) = f^T(4 + 4) = f^T(8) = \sup\{m; m \in N \text{ and } m \leq \sqrt{8}\} = 2_\mu$ because $f^T(3) = 9 > 8$ or $3 > \sqrt{8}$.

$2_\mu + 3_\mu = f^T(f_T(2_\mu) + f_T(3_\mu)) = f^T(4 + 9) = f^T(13) = \sup\{m; m \in N \text{ and } m \leq \sqrt{13}\} = 3_\mu$ because $f^T(4) = 16 > 13$ or $4 > \sqrt{13}$.

$10_\mu + 11_\mu = f^T(f_T(10_\mu) + f_T(11_\mu)) = f^T(100 + 121) = f^T(221) = \sup\{m; m \in N \text{ and } m \leq \sqrt{221}\} = 14_\mu$ because $f^T(15) = 225 > 221$ and $f^T(14) = 196 < 221$.

In a similar way, we find that:

$2_\mu + 11_\mu = 11_\mu$, $3_\mu + 11_\mu = 11_\mu$, $4_\mu + 11_\mu = 11_\mu$, but $5_\mu + 11_\mu = 12_\mu$, $6_\mu + 11_\mu = 12_\mu$, $7_\mu + 11_\mu = 13_\mu$, $8_\mu + 11_\mu = 15_\mu$, and $11_\mu + 11_\mu = 15_\mu$.

$2_\mu \circ 2_\mu = f^T(f_T(2_\mu) \cdot f_T(2_\mu)) = f^T(4 \cdot 4) = f^T(16) = 4_\mu$.

$2_\mu \circ 3_\mu = f^T(f_T(2_\mu) \cdot f_T(3_\mu)) = f^T(4 \cdot 9) = f^T(36) = 6_\mu$.

In general, we have $n^2 \cdot m^2 = (n \cdot m)^2$. Consequently, in this arithmetic $A$, multiplication of numbers is the same as in the Diophantine arithmetic, i.e., $n_\mu \circ m_\mu = (n \cdot m)_\mu$, while addition is essentially different. As we have seen, two times two is still four, while two plus two is only two.

**Example 3.** Let us take as a generator for a projective arithmetic $A = (N; +, \circ, \leq)$ a simple function, such as $10x$.

$2_\mu + 2_\mu = f^T(f_T(2_\mu) + f_T(2_\mu)) = f^T(20 + 20) = f^T(40) = 4_\mu$.

$2_\mu + 3_\mu = f^T(f_T(2_\mu) + f_T(3_\mu)) = f^T(20 + 30) = f^T(50) = 5_\mu$.

This is a general case for addition of numbers in $A$ as $n_\mu + m_\mu = f^T(f_T(n_\mu) + f_T(m_\mu)) = f^T(10n + 10m) = f^T(10(n + m)) = (n + m)_\mu$.

At the same time, we have:

$2_\mu \circ 2_\mu = f^T(f_T(2_\mu) \cdot f_T(2_\mu)) = f^T(20 \cdot 20) = f^T(400) = 40_\mu$.

$2_\mu \circ 3_\mu = f^T (f_T (2_\mu) \cdot f_T (3_\mu)) = f^T (20 \cdot 30) = f^T (600) = 60_\mu$.

This is a general case for addition of numbers in $A$ as $n_\mu \circ m_\mu = f^T (f_T (n_\mu) \cdot f_T (m_\mu)) = f^T (10n \cdot 10m) = f^T (10(n \cdot m)) = (n \cdot m)_\mu$.

In this arithmetic $A$, addition of numbers is the same as in the Diophantine arithmetic, while multiplication is essentially different. As we have seen, two plus two is still four, while two times two is equal to forty.

**Theorem 3.** Both operations, addition + and multiplication ∘, are commutative in any projective arithmetic.

**Remark 2.** In some projective arithmetics, addition + or/and multiplication ∘ can be non-associative.

There are other operations in $A$. For instance, we have two $n$-ary operations:

1. $\Sigma^n (a_{\mu,1}, a_{\mu,2}, \ldots, a_{\mu,n-1}, a_{\mu,n}) = f^T (f_T (a_{\mu,1}) + f_T (a_{\mu,2}) + \ldots + f_T (a_{\mu,n}))$;

2. $\Pi^n (a_{\mu,1}, a_{\mu,2}, \ldots, a_{\mu,n-1}, a_{\mu,n}) = f^T (f_T (a_{\mu,1}) \cdot f_T (a_{\mu,2}) \cdot \ldots \cdot f_T (a_{\mu,n}))$.

**Example 4.** Let us take a projective arithmetic $A = (N; +, \circ, \leq)$ with the generator $f(x) = x^2$, and look how $n$-ary operations are performed.

$\Sigma^n (2_\mu, 2_\mu, 2_\mu) = f^T (4 + 4 + 4) = f^T (12) = 3_\mu$.

$\Sigma^n (2_\mu, 2_\mu, 3_\mu) = f^T (4 + 4 + 9) = f^T (17) = 4_\mu$, but

$\Sigma^n (1_\mu, 1_\mu, 1_\mu, 1_\mu, 1_\mu, 3_\mu) = f^T (1 + 1 + 1 + 1 + 1 + 9) = f^T (14) = 3_\mu$.

$\Pi^n (2_\mu, 5_\mu, 8_\mu) = f^T (4 \cdot 25 \cdot 64) = f^T (6400) = 80_\mu = (2 \cdot 5 \cdot 8)_\mu$.

This is a general case for the operation $\Pi^n$ as we have the equality

$$m_1^{\,2} \cdot m_1^{\,2} \cdot \ldots \cdot m_n^{\,2} = (m_1 \cdot m_1 \cdot \ldots \cdot m_n)^2$$

In the Diophantine arithmetic, we have $\Sigma^n (a_1, a_2, \ldots, a_{n-1}, a_n) = a_1 + a_2 + \ldots + a_n$. However, in some projective arithmetics, addition is not associative, and we do not have the corresponding identity $\Sigma^n (a_{\mu,1}, a_{\mu,2}, \ldots, a_{\mu,n-1}, a_{\mu,n}) = a_{\mu,1} +$

$a_{μ, 2} + \ldots + a_{μ, n}$. Thus, it is an important question in what projective arithmetics addition is associative.

**Theorem 4.** Addition + in a projective arithmetic $A = (M ; +, \circ, \leq)$ is associative if and only if the function $f$ is piecewise linear.

In addition to $\leq$ and $<$, there are other kinds of order relations in $A$:

1) $a_μ \ll b_μ$ means that $a_μ$ is *much less* than $b_μ$ and in this case, $b_μ$ is *much larger* than $a_μ$ :

$$a_μ \ll b_μ \text{ if and only if } b_μ + a_μ = b_μ.$$

2) $a_{μ,1}, a_{μ,2}, \ldots, a_{μ,n-1} \ll_n b_μ$ means that the group $a_{μ,1}, a_{μ,2}, \ldots, a_{μ,n-1}$ is *much less* than $b_μ$ and in this case, $b_μ$ is *much larger* than the group $a_{μ,1}, a_{μ,2}, \ldots, a_{μ,n-1}$:

$$a_{μ,1}, a_{μ,2}, \ldots, a_{μ,n-1} \ll_n b_μ \text{ if and only if } \Sigma^n (a_{μ,1}, a_{μ,2}, \ldots, a_{μ,n-1}, b_μ) = b_μ.$$

3) $a_μ \lll b_μ$ means that $a_μ$ is *much much less* than $b_μ$ and in this case, $b_μ$ is *much much larger* than $a_μ$ :

$$a_μ \lll b_μ \text{ if and only if } b_μ \circ a_μ = b_μ.$$

4) $a_{μ,1}, a_{μ,2}, \ldots, a_{μ,n-1} \lll_n b_μ$ means that the group $a_{μ,1}, a_{μ,2}, \ldots, a_{μ,n-1}$ is *much much less* than $b_μ$ and in this case, $b_μ$ is *much much larger* than the group $a_{μ,1}, a_{μ,2}, \ldots, a_{μ,n-1}$:

$$a_{μ,1}, a_{μ,2}, \ldots, a_{μ,n-1} \lll_n b_μ \text{ if and only if } \prod^n (a_{μ,1}, a_{μ,2}, \ldots, a_{μ,n-1}, b_μ) = b_μ.$$

It is possible that $\ll$ is a total relation on $M$.

**Example 5.** Let us take $f_T(n) = 2^{2^n}$ as a generator for a projective arithmetic. Then for any $n$, we have

$$f_T(n+1) + f_T(n) = 2^{2^{n+1}} + 2^{2^n} < 2^{2^{n+1}} + 2^{2^{n+1}} < 2^{2^{n+1}} \cdot 2^{2^{n+1}} = 2^{2^{n+1} + 2^{n+1}} = 2^{2 \cdot 2^{n+1}} = 2^{2^{n+2}} = f_T(n+2)$$

Consequently, we have the following relations $1_\mu \ll 2_\mu \ll 3_\mu \ll \ldots \ll n_\mu \ll (n+1)_\mu \ll \ldots$ in the projective arithmetic with the projector $f_T(n)$.

**Theorem 5.** Relation $\ll$ in a projective arithmetic $A = (M; +, \circ, \leq)$ is transitive and disjunctively asymmetric, i.e., only one relation $xQy$ or $yQx$ is valid for all different elements $x$ and $y$ from $M$.

The proof is based on the following lemma.

**Lemma 1.** For any number $a$ from a projective arithmetic $A$, we have $a \ll Sa$ if $f_T(Sa) + f_T(a) < f_T(SSa)$.

**Definition 4.** a) A binary relation P on a set $X$ is called *compatible from the right (from the left)* with a binary relation Q on $X$ if $P \circ Q \subseteq P$ ($Q \circ P \subseteq P$).

b) A binary relation P on a set $X$ is called *compatible* with a binary relation Q on $X$ if P is compatible both from the right and from the left with Q.

**Theorem 6.** Relation $\ll$ in a projective arithmetic $A = (M; +, \circ, \leq)$ is compatible with the order $\leq$ in $M$.

**Remark 3.** Results of Theorems 5 and 6 are not always true for the relation $\lll$ as the following example demonstrates.

**Example 4.** Let us take as a generator for a projective arithmetic the following function:

$$f_T(n) = \begin{cases} 2^{2^n} & \text{for } n = 1, 2, 3, 4 ; \\ 2^{n+12} & \text{when } n > 4 \end{cases}$$

Then we have

$$2_\mu \circ 3_\mu = f^T(f_T(2) \cdot f_T(3)) = f^T(2^{2^2} \cdot 2^{2^3}) = f^T(2^4 \cdot 2^8) = f^T(2^{12}) < 4_\mu = f^T(2^{2^4}) = f^T(2^{16})$$

.

Thus, $2_\mu \circ 3_\mu = 3_\mu$ and $2_\mu <<< 3_\mu$.

At the same time, we have $2_\mu \circ 5_\mu = f^T(f_T(2) \cdot f_T(3)) = f^T(2^{2^2} \cdot 2^{17}) = f^T(2^4 \cdot 2^{17}) = f^T(2^{21}) = 9_\mu$. Consequently, it is not true that $2_\mu <<< 5_\mu$ although $2_\mu < 3_\mu < 5_\mu$. It means that the relation $<<<$ is not compatible from the right with the relation $<$.

However, assuming that the generator $f(x)$ is a monotonous function, we have the following result.

**Theorem 7.** Relation $<<<$ in a projective arithmetic $A = (M; +, \circ, \leq)$ is compatible from the left with the order $\leq$ in $M$.

Let $A = (M; +, \circ)$ be a projective prearithmetic, for which the second condition from the definition 4 is true and $f_T(1_\mu) = 1$.

**Theorem 8.** For any $n \in N$ and any $a_\mu, b_\mu \in A$, from $1_\mu, 1_\mu, \ldots, 1_\mu <<_n a_\mu$ follows $1_\mu, 1_\mu, \ldots, 1_\mu <<_n b_\mu$ if and only if for any elements $a$ and $b$ from $M$ from $a \leq b$ it follows $f_T(Sa) - f_T(a) \leq f_T(Sb) - f_T(b)$.

## 5. Conclusion

In the paper, explicit constructions for prearithmetics and non-Diophantine arithmetics are gives. Properties of such arithmetics are studied. Such properties can be related to the basic features of nature. Some physicists (cf., for example, (Zeldovich, *at al*, 1990)) emphasized that fundamental problems of modern physics are dependent on our ways of counting. This idea correlates with problems of modern physical theories in which physical systems are described by chaotic processes. Taking into account the fact that chaotic solutions are obtained by computations, physicists ask (Cartwrite and Piro, 1992; Gontar, 1997) whether chaotic solutions of the differential equations, which model different physical systems, reflect the dynamic laws of nature represented by these equations or whether they are solely the result of an extreme sensitivity of these solutions to numerical procedures and computational errors.

It is even clearer that properties of non-Diophantine arithmetics, which reflect the way people count, influence functioning of economy and are important for economical models (cf., for example, (Tolpygo, 1997)). Thus, it would be useful to build models of economical systems and processes using not on the Diophantine arithmetic but an appropriate non-Diophantine arithmetic.

One of the interesting properties of projective arithmetics is that they allow us to formalize and make rigorous such concepts as *much smaller* (<<) and *much larger* (>>), which are frequently used in theoretical physics and other areas of applied mathematics on the intuitive level without a proper mathematical formalization. Namely, we have the following definition:

**A number *m* is *much smaller* than a number *n* ( *m* << *n*) if *n* + *m* = *n*.**

In this case, the number *n* is *much larger* than the number *m* ( *n* >> *m*)

Non-Diophantine arithmetics also allow one to eliminate several inconsistencies and misconceptions related to arithmetic. For instance, Rosinger (2008) explains that "we have been doing inconsistent mathematics for more than half a century by now, and in fact, have quite heavily and essentially depended on it in our everyday life. Indeed, electronic digital computers, when considered operating on integers, which is but a part of their operations, act according to the system of axioms given by

- (PA) : the usual Peano Axioms for *N*,

plus the ad-hock axiom, according to which

- (MI) : there exists *M* in *N*, M >> 1, such that M + 1 = M

Such a number *M*, called "machine infinity", is usually larger than 10100, however, it is inevitably inherent in every electronic digital computer, due to obvious unavoidable physical limitations. And clearly, the above mix of (PA) + (MI) axioms is inconsistent. Yet we do not mind flying on planes designed and built with the use of such electronic digital computers."

In a similar way, Meyer and Mortensen (1984) built various inconsistent models of arithmetic, while Priest (1997; 2000) developed axiomatic systems for inconsistent arithmetics.

Even before Priest, Van Bendegem (1994) developed an inconsistent axiomatic arithmetic by changing the Peano axioms so that a number that is the successor of itself exists. The fifth Peano axiom states that if $x = y$ then $x$ and $y$ are the same number. In the system of Van Bendegem, starting from some number $n$, all its successors will be equal to $n$. Then the statement $n = n + 1$ is considered as both true and false at the same time. This makes the new arithmetic inconsistent.

There are two basic ways to deal with inconsistencies: one is to elaborate an inconsistent system and try to work with it and another way is create new mathematical structures, eliminating inconsistencies. The book "Mathematicians Think" of William Byers (2007) has the subtitle "Using Ambiguity, Contradiction, and Paradox to Create Mathematics" because mathematical reasoning, according to Byers, is not completely algorithmic, computational or based on proof systems, but primarily uses creative ideas to shed new light on mathematical objects and structures, propelling in such a way mathematical progress. The central role in the emergence of creative ideas ambiguities, contradictions, and paradoxes play. According to Byers, one of the main kinds of contradictions is existence of two seemingly contradictory perspectives in a mathematical problem. For instance, Peano axioms imply infiniteness of the arithmetic, while the existence of the largest number implies it finiteness. However, this paradox vanishes with the discovery of non-Diophantine arithmetics. Actually all these inconsistencies and contradictions exist only in the absence of non-Diophantine arithmetics. For instance, the machine arithmetic analyzed by Rosinger (2008) does not satisfy

Peano axioms because it is a non-Diophantine arithmetic, while peano axioms formalize the Diophantine arithmetic. That is why the machine arithmetic satisfies the axioms of the corresponding non-Diophantine arithmetic and electronic digital computers, when considered operating on integers, which is a part of their operations, act according to this system of axioms.

The rule $n = n + 1$ of Van Bendegem (1994) is natural for many non-Diophantine arithmetics and causes no inconsistencies and contradictions there, meaning simply $1 \ll n$.

It is possible to compare this situation with artificially derived inconsistencies and contradictions with the time when people knew only natural numbers and positive fractions. Getting information about negative numbers, mathematicians who lived at that time would be able to build an inconsistent formal system by taking two "axioms":

- Only positive numbers exist;

and

- There are negative numbers.

Naturally these "axioms" give a contradiction. Now we know that the first "axiom" is valid only for natural numbers and positive fractions if we consider numbers known at that time. Thus, integer numbers combine both positive and negative numbers without any inconsistency.

Mathematicians who lived in the 19$^{th}$ century and earlier were able also to build an inconsistent formal system in geometry, combining together two sets of axioms:

- All postulates of the Euclidean geometry;

and the postulate which is true for the geometry on a sphere, which is considered the geometrical model of the Earth

- Any two straight lines intersect with one another.

Now we know that spherical geometry is non-Euclidean and does not have any contradictions in it.